\documentclass[12pt]{article}
\usepackage{amsthm}
\usepackage[T1]{fontenc}
 \usepackage{graphicx}
\usepackage{color}

\usepackage{mathtools}
\usepackage{amsmath,amssymb,amsfonts}        
\usepackage[mathscr]{euscript}
\usepackage{dsfont}

\newtheorem{definition}{Definition}[section]
\newtheorem{theorem}{Theorem}[section]

\newtheorem{comment}{Comment}[section]

\newcommand{\R}{\mathds{R}}

\begin{document}




\title{Solving multi-resource allocation and location problems in disaster management through linear programming}

\maketitle
\author{B. Bodaghi, Faculty of Science, Engineering and Technology,\\ Swinburne University of Technology, PO Box~218, Hawthorn, Victoria,\\ \textsc{Australia}, bbodaghi@swin.edu.au}

\author{N. Sukhorukova, {Faculty of Science, Engineering and Technology,\\ Swinburne University of Technology, PO Box~218, Hawthorn, Victoria,\\ \textsc{Australia} and Centre for Informatics and Applied Optimization,\\ Federation University Australia, nsukhorukova@swin.edu.au}

\begin{abstract}
In this paper we propose a new efficient linear programming based approach for multi-resource  allocation and location problems in  disaster management. Such problems require an integer solution and therefore, in most cases, the computations rely on integer and mixed-integer linear programming solvers. In general, these solvers can not handle large scaled problem. 
In this paper we demonstrate that there exists a large class of disaster management problems whose exact solutions can be obtained by applying the simplex method (linear programming). 
The results of numerical experiments are provided.
Another important contribution of this paper is related to general cluster analysis and allocation. Namely, we demonstrate that the classical $k$-medoid clustering method can be implemented using  linear programming techniques (simplex method) without relying on integer solvers.   
\end{abstract}

{\bf Keyword:} 
linear programming,  integer programming, optimisation, simplex method, disaster management, resource allocation, relief centre location


MSC codes: 90C06,  90C08,  90C11,  90C05


\section{Introduction}
\label{sec:intro}

Many disaster management approaches are based on optimisation techniques. In particular, the problems are formulated as mathematical programs and then solved using available optimisation techniques {\cite{ZHOU2018567}}. In the simplest case,  problems can be formulated as  linear or convex  programs, while some others require more   advanced techniques.

A large class of applications can be formulated as mixed-integer linear programs. In these problems, all or some of the variables are integers. This makes the corresponding optimisation problems much harder to solve. Moreover, in some cases, we can not even guaranty that an obtained solution is optimal.

One approach to solve these problems is to apply standard linear optimisation methods  and then round the solution to the nearest integer. This approach is not always very precise, especially when the integer variables are binary (that is, can only take values from $\{0,1\}$). In this paper we derive a rich class mixed-integer linear problems whose optimal solution can be obtained by applying the simplex method, a powerful methods developed for classical linear programs. 

There are several advantages of our approach. First of all, we can guaranty that the obtained solution is integer and optimal. Second, the simplex method is part of most linear programming packages, where  mixed-integer programming solvers may not be available. Third,  linear programming techniques can deal with very large problems,   while the applicability of mixed-integer techniques to same size problems may not be possible. 

In this paper we consider three applications taken from the field of disaster management. There are many more applications where our approach is applicable. The list of such applications is rich and goes well beyond the area of disaster management. 

The paper is organised as follows. In section~\ref{sec:disastermanagement} we give an overview of disaster management problems. Then, in section~\ref{sec:LinearProg}, we formulate the most important results from the area of linear programming. In section~\ref{sec:SimplexUnimodular} we discuss the application of the classical simplex method to integer programming problems whose constraint matrices are totally unimodular. In section~\ref{sec:modelling} we formulate our main results. In section~\ref{sec:ne} we provide the results of the numerical experiments. Finally, in section~\ref{sec:conclusions} we provide our conclusions and discussions.

\section{Disaster management }\label{sec:disastermanagement}
%
%
Disaster management contains five general phases, namely: prediction, warning, emergency relief, rehabilitation, and reconstruction \cite{lin2006integrated}. The related activities which are conducted in these phases usually contain mitigation and preparedness, response, and recovery \cite{jayaraman1997managing}.

In emergency relief operation, resource allocation problem (RAP) is often a complex challenge due to a number of issues, such as dealing with crucial demands, time criticality, competing priorities, the extent of availabilities, and different constraints and uncertainties \cite{sheu2007challenges}. Emergency resources can be grouped  into non-expendable and expendable resources~\cite{Bodaghi2018,SHAHPARVARI2018886}. Non-expendable resources are non-consumable and renewable and might include emergency personnel and volunteers. Expendable resources are consumable and cannot be renewed in the emergency, or the recovery cost surpasses the original value of the resource (for instance, medical supplies, water, food and fuel).  Failure to assign adequate resources in a timely manner has been the main cause of adverse impacts in disaster situations \cite{rolland2010decision,lei2015personnel,najafi2013multi}. Efficient reduction of losses and fatalities in disaster-struck locations are significantly dependent on the rapid deployment of resources for emergency relief operations.

A RAP models consider assignment of resources to task from relief centres without determining the sequence or flow of resources along arcs \cite{caunhye2012optimization}. A RAP problem in disaster management is usually formulated as a mixed-integer linear programming (MILP) problem~\cite{Bodaghi2018,caunhye2012optimization,fiedrich2000optimized}. Most of the methodologies in natural disaster response phase focus on developing a MILP model just for one type of the resources (or just integer linear programming, ILP, when all the variables are integers). A literature review reveals  that RLP models with different kinds of resources  in disaster management have been discussed quite rarely (\cite{lei2015personnel,lee2013operations,Bodaghi2018}. 

In addition, in most cases it is assumed that one relief centre may supply several demand points, while every demand point can receive supplies from just one relief centre. There are two main advantages of this approach. First of all, this rule makes decision maker\rq{}s tasks easier to implement. Second, this approach may minimise the number of vehicles involved. The main disadvantage of this model is that the total costs are higher than they may be in the case when several relief centres are allowed to supply the same demand point. Later in this paper we  demonstrate that the removal of this assumption leads to MILP problems that can be solved using standard linear programming tools, that are very advanced and efficient, while the obtained solution is optimal for the original MILP problem.

{Furthermore, finding the location of relief centres (e.g. medical centres, distribution centres) are vital in disaster management~\cite{BOONMEE2017485}. This problem can be addressed as a facility location problem (e.g.~\cite{Balcik}) or a clustering problem (e.g.~\cite{Paul2012,Bodaghi2017}). To find a solution for the location of relief centres normally a MILP model is formulated for small cases, and using several heuristics approaches for solving the larger cases. Heuristic algorithms are producing a non-optimal solution and the result of them are poor quality when benchmarked with exact approach (e.g. Linear Programming) \cite{BOONMEE2017485}}.

\section{Linear Programming and Transportation problem}\label{sec:LinearProg}

\subsection{Linear Programming}\label{subsec:LinearProg}

A general {linear programming problem} or {LPP} can be formulated as follows:
\begin{equation}\label{eq:LPP_obj}
\text{minimise}~{\bf c}^T{\bf x}
\end{equation}

\begin{align}
\text{subject~to~}& {\bf a}_i^T{\bf x}\geq b_i,&i\in M_1;\label{eq:LPP_constraints1}\\
                & {\bf a}_i^T{\bf x}\leq b_i,&i\in M_2;\label{eq:LPP_constraints2}\\
                &{\bf a}_i^T{\bf x}= b_i,&i\in M_3;\label{eq:LPP_constraints3}\\
                &x_j\geq 0, &j\in N_1;\label{eq:LPP_constraints4}\\
                &x_j\leq 0, &j\in N_2\label{eq:LPP_constraints5},
\end{align}
where
\begin{itemize}
\item  $x_1,\dots,x_n$ are the {decision variables};
\item ${\bf c}=(c_1\dots,c_n)$ is the {cost vector};
\item function ${\bf c}^T{\bf x}$ in~(\ref{eq:LPP_obj}) is the {objective (cost) function};
\item  the set of equalities and inequalities~(\ref{eq:LPP_constraints1})-(\ref{eq:LPP_constraints5}) are  the {constraints}.
\end{itemize}

\begin{definition}
A vector ${\bf x}$ satisfying all the constraints is called a {feasible solution}. The set of all feasible solutions is called the {feasible set}.
\end{definition}
\begin{definition}
 A feasible solution ${\bf x}^*$ that minimises the objective function is called  an {optimal feasible solution} or just an {optimal solution}.
 \end{definition}
 \begin{definition}
 A liner programming problem is called feasible if there exists at least one feasible solution. Otherwise, the problem is infeasible.
 \end{definition}

There are many efficient methods for solving linear programming problems: the simplex method, originally developed by Dantzig in 1947~\cite{GB}, interior point methods~\cite{KN,JS}. It has been demonstrated in 1972 by~\cite{KM} that  the worst-case complexity of the simplex method as formulated by Dantzig is exponential time. The  worst-case complexity of interior point methods is polynomial. Despite these results, the simplex method is remarkably efficient and included in most linear programming packages.  

In many practical problems it is beneficial (where it is possible) to reformulate an integer or mixed-integer linear programming problem as an LPP and solve it then using linear programming techniques. This normally involves a significant increase in the number of decision variables and/or constraints.

\subsection{Transportation Problems}\label{subsect:Transportation}

Transportation problems form a special class of LPPs. In most operations research textbooks, this problem is formulated in a way related to transportation. It is possible, however, to apply the same type of approach to other types of problems, including location analysis and allocation problems. An excellent overview of such problems and also more advanced models  can be found in~\cite{ZviDrezner,Location:2015}.  

Consider an example of transportation problem formulation. Goods are produced at $m$ factories (also called sources)
$$S_1,\dots, S_m$$ 
and  sold at $n$ markets (also called destinations):
$$D_1, \dots, D_n.$$ 
The supply available at source $S_i$ is $s_i\geq 0$ units, the demand at destination~$D_j$ is $d_j\geq 0$ units  and the transportation cost of one unit from
$S_i$ to $D_j$ is $c_{ij}\geq 0$. We have to identify which sources should supply which destinations to minimise total  transportation costs.

Let $x_{ij}$ be the number of units to be sent from $S_i$ to $D_j$. Then the corresponding optimisation problem can be formulated as follows:
\begin{equation}
\min \sum_{i=1}^m\sum_{j=1}^nc_{ij}x_{ij}
\end{equation}
subject to 
\begin{equation}
\sum_{j=1}^n x_{ij}\leq s_i,~i=1,\dots,m;\\
\sum_{i=1}^m x_{ij}\geq d_j,~j=1,\dots,n;\\ 
\end{equation}
\begin{equation}
x_{ij}\geq 0,~i,j=1,\dots,n.
\end{equation}

It is well known that if all supplies  and demands are integers, then there exists an optimal solution $x_{ij}$, which is integer. This is important for many applications where the units (for example, computers, cars, people) can not be split. In general,  integer and mixed-integer programming problems are much harder to solve than linear programming problems. It is also well known that the simplex method  applied to a transportation problem, terminates at an optimal solution that is also integer. 

There are other types of LPPs whose optimal solutions, reached by applying the simplex method, are integers. 
 In the next section we provide a brief description of the simplex method and identify a broad class of linear problems where an integer optimal solution  can be obtained by applying standard linear programming techniques. 

\section{The simplex method and its application to integer problems whose constraint matrices are totally unimodular}\label{sec:SimplexUnimodular}

The simplex method is a powerful linear programming algorithms.    The algorithm starts at a feasible vertex of the constraint polyhedron.  Then it move to an adjacent vertex, where the objective function value is at least as good as it is at the original vertex. Since the number of vertices  is finite, this method will terminate in a finite number of steps. If all the vertices of the constrain polyhedron have integer coordinates, then an optimal solution (can be more than one) is integer. 

Consider a system of linear inequalities
\begin{equation}\label{eq:LPforSimplex}
{\bf Ax}\leq{\bf b},
\end{equation}
where ${\bf A}$ is an $m\times n$ matrix, ${\bf x}\in \R^n$ and ${\bf b}\in \R^m$. It is well-known that all possible solutions to a system of linear inequalities form a polyhedron. This polyhedron is an empty set if the system does not a solution.  

A comprehensive overview on integer programming can be found in~\cite{Schrijver:1986}. Theorem~19.3 of this book (p. 268) covers the conditions when the vertices of the feasible sets are integers and therefore an optimal solution found at a vertex is integer.  In this paper, we only use conditions (i), (iii) and (iv) of the  theorem (originally proved in~\cite{HoffmanKruskal:1956} and \cite{Ghouila:1962}). A simplified version of this theorem, formulated for this study, is as follows.
\begin{theorem}\label{thm:main}
Let $A$ be a matrix with entries $\{0,1,-1\}$. Then the following are equivalent:
\begin{enumerate}
\item matrix~$A$ is totally unimodular, that is each square submatrix of~$A$ has determinant $0$, $1$ or $-1$; 
\item for all integral vectors~$a$, $b$, $c$ and $d$ the polyhedron 
$$\{x| c\leq x\leq d,~a\leq Ax\leq b\}$$
has only integral vertices;
\item each collection of columns of $A$ can be split into two parts so that the sum of the columns in one part minus the sum of the columns in the other part is a vector with entries~$\{0, 1, -1\}$;
\end{enumerate}
\end{theorem}

The first condition of Theorem~\ref{thm:main} is usually used as a definition for totally unimodular matrices. 
The class of totally unimodular matrices is closed under a number of operations. We need the following ones:
\begin{itemize}
\item
  transposition;
  \item multiplication a row (column) by $-1$.
 \end{itemize}
Also, matrix $A$ is totally unimodular if and only if matrix~$[I A]$ (where $I$ is an identity matrix of the corresponding dimension) is totally unimodular.

In the next section we formulate three types of integer problems appearing in disaster management. Then we demonstrate that the corresponding system matrices are totally unimodular and therefore their optimal solutions (reached at vertices) are integers  and hence   can be obtained by applying the simplex method.

\section{Mathematical Modelling}\label{sec:modelling}
\subsection{{Expendable} resources}\label{subsec:consumable}

The problem of allocating expendable resources can be formulated as a Transportation Problems. It is enough to think about incident points as ``Markets\rq{}\rq{} (each market demand corresponds to the corresponding incident point demand), while the relief centres are ``Factories\rq{}\rq{} (each factory capacity corresponds to the processing centre capacity). The transportation costs are ``processing and transportation time\rq{}\rq{}. 

The feasible set of this problem is as follows (without sign constraints and integer requirement):
\begin{equation}\label{eq:exp}
\left[
\begin{matrix}
-I_n&-I_n&-I_n&\dots &-I_n\\
e_n&0_n&0_n&\dots &0_n\\
0_n&e_n&0_n&\dots &0_n\\
0_n&0_n&e_n&\dots &0_n\\
\vdots&\vdots&\vdots&\ddots &\vdots\\
0_n&0_n&0_n&\dots &e_n\\
\end{matrix}
\right]X\leq b,
\end{equation}
where 
\begin{itemize}
\item $b\in\R^{(n+m)}$ represents the corresponding demands and supplies and therefore $b$ is integral;
\item $X\in\R^{mn}$ is the vector of decision variables;
\item $I_n$is an identity matrix of size $m$;
\item $e_n(1,1,\dots,1)\in\R^n$;
\item $0_n\in\R^n(0,0,\dots, 0)$;
\item the system matrix $A\in\R^{(n+m)\times(mn)}$. 
\end{itemize}

\begin{theorem}\label{thm:trans_unimod}
The system matrix $A$ from (\ref{eq:exp}) is totally unimodular.
\end{theorem}
{\bf Proof:} Consider matrix ~$B$ obtained from $A^T$ by multiplying the first $n$ columns of $A^T$ by $-1$. Matrix~$A$ is totally modular if and only if matrix~$B$ is totally unimodular.

   Assign the first $n$   columns of $B$ to part~I and the remaining columns to part~II and assume that one or more columns may be removed from the total collection of columns. The sum of the columns in part~I is an $(mn)$-dimensional vector~$S_1$ whose components are~$0$ or $1$. The sum of the columns in  part~I is an $(mn)$-dimensional vector~$S_2$ whose components are~$0$ or $1$. Therefore the components of~$S_1-S_2$ are $0$, $1$ or $-1$ and hence, by Theorem~\ref{thm:main}, we conclude that matrices $B$ and $A$ are totally unimodular and all the vertices of the feasible set have integer coordinates.
   
   \hskip300pt $\square$
   
   We can also conclude that  an optimal integer solution to this problem can be found by applying the simplex method.

%
%

\subsection{Non-expendable resources}\label{subsec:nonconsumable}

In the case of non-expendable resources, the problem can also be formulated as an integer programming problem, where some of the summation constraints from a classical transportation problem  are replaced with maximisation. This problem is not a classical transportation problem, but it can be formulated as an LPP. It can be demonstrated that the applications of the simplex method also leads to an integer optimal solution.

A mixed-integer formulation for the case of non-consumable resources is as follows
\begin{equation}\label{eq:trans_obj_c}
\min \sum_{i=1}^m\sum_{j=1}^n c_{ij}x_{ij}
\end{equation}
subject to
\begin{equation}\label{eq:trans_con_1_c}
\sum_{i=1}^{m}x_{ij}\geq d_j,~j = 1,\dots, n;
\end{equation}
\begin{equation}\label{eq:trans_con_2_c}
\max_{j=1,\dots,n}x_{ij} \leq s_i,~i = 1,\dots,m;
\end{equation}
\begin{equation}\label{eq:trans_con_box_c}
x_{ij}\geq 0,~i=1,\dots,n,~j=1,\dots,m,
\end{equation}
\begin{equation}\label{eq:trans_int}
x_{ij}~{\rm are~integers},~i=1,\dots,m,~j=1,\dots,n,
\end{equation}
where $d_i,~i=1,\dots,m$ are incident point demands and $s_j,~j=1,\dots,n$ are  relief centre capacities.
A relaxation of this problem, obtained by removing the last constraint~(\ref{eq:trans_int}),  can be formulated as an LPP by replacing constraints~(\ref{eq:trans_con_2_c}) with equivalent  systems of linear inequalities:
 \begin{equation}\label{eq:trans_con_2_cc}
x_{ij} \leq s_j,~j = 1,\dots,n,~i=1,\dots,m.
\end{equation}

The feasible set of this problem can be formulated as follows (without sign constraints and integer requirement):
\begin{equation}\label{eq:non-exp}
\left[
\begin{matrix}
-I_n&-I_n&-I_n&\dots &-I_n\\
      &        & I_{mn} &   &  \\
\end{matrix}
\right]X\leq b,
\end{equation}
where 
\begin{itemize}
\item $b\in\R^{(n(m+1))}$ represents the corresponding demands and supplies and therefore $b$ is integral;
\item $X\in\R^{mn}$ is the vector of decision variables;
\item $I_n$ is an identity matrix of size $n$;
\item $I_{mn}$ is an identity matrix of size $mn$;
\item the system matrix $A\in\R^{n(m+1)\times(mn)}$. 
\end{itemize}

\begin{theorem}
The system matrix $A$ from~(\ref{eq:non-exp}) is totally unimodular.
\end{theorem}
{\bf Proof:}
 $A$ is totally unimodular if and only if matrix
  \begin{equation}
B=\left[
\begin{matrix}
I_n&I_n&I_n&\dots &I_n\\
\end{matrix}
\right]
\end{equation}
is totally unimodular. Matrix $B$ is totally unimodular if and only if $I_n$ is totally unimodular. Indeed, one can assign any collection of columns of $I_n$ to part~I and the remaining columns to part~II. The difference of the corresponding columns sums contains $1$ and $-1$ as the components. 
%
   
   \hskip300pt $\square$

  Therefore, similar to the case of expendable resources, we can reduce an integer linear programming problem to an LPP whose vertices are integral.

\subsection{Relief centre location}


In this application, the distance matrix between all the incident points is given. The goal is to select $k$ points in such a way that, after assigning all the remaining points to the nearest selected point (cluster centre), the total sum of distances between the points and centres is minimal. Each cluster centre is a relief centre, whose optimal location (selection among the incident points) is the objective. In this application, we assume that the demand of the incident points can be covered regardless of the allocation, since the main objective is to minimise the total distance. This kind of clustering is called~$k$-medoid, was first proposed in~\cite{KaufmanRousseeuw}.

Assume that there are $n$ demand points in total and the distance matrix $${\bf D}=\{d_{ij}\},~i=1,\dots,n,~j=1,\dots,n.$$ 
It is easy to see that this matrix is symmetric and its main diagonal consists of zeros. The goal is to select~$k$ points as relief centres. The decision variables are binary: 
$$x_{ij}\in\{0,1\},~i=1,\dots,n,~j=1,\dots,n$$ and $y_i\in\{0,1\},~i=1,\dots,n.$
Variable $y_i$ is $1$ if incident  point~$i$ is treated as a relief centre, otherwise, this variable is zero. Variable $x_{ij}$ is $1$ if incident point~$i$ was assigned to point~$j$.
The corresponding optimisation problem is as follows:
\begin{equation}\label{eq:kmedoid_obj}
\min \sum_{i=1}^n\sum_{j=1}^n d_{ij}x_{ij}
\end{equation}
subject to
\begin{equation}\label{eq:kmedoid_c1}
\sum_{i=1}^{n}x_{ij}=1,~j = 1,\dots, n;
\end{equation}
 \begin{equation}\label{eq:kmedoid_c2}
x_{ij}\leq y_i,~i,j=1,\dots, n;
\end{equation} 
 \begin{equation}\label{eq:kmedoid_c3}
\sum_{i=1}^{n}y_{i}=k;
\end{equation} 
 \begin{equation}\label{eq:kmedoid_c4}
x_{ij},~y_i\in\{0,1\},~i,j=1,\dots,n.
\end{equation} 
Constraints~(\ref{eq:kmedoid_c1}) ensure that each incident point is assigned to a single  relief centre.  Constraints~(\ref{eq:kmedoid_c2}) ensure that an incident point~$i$ can only be assigned to an incident point~$j$ if this point is also a  relief centre. Finally, constraint~(\ref{eq:kmedoid_c3}) ensures that exactly~$k$ points are selected as  relief centres.
It is clear that problem~(\ref{eq:kmedoid_obj})-(\ref{eq:kmedoid_c4}) is an integer programming problem and, in general, it is not easy to solve this problem.

\begin{theorem}\label{thm:kmedoid}
All the vertices of the feasible set in $k$-medoid method have integer coordinates.
\end{theorem}
{\bf Proof:}
 The set of constraints contains $(n+1)$~equalities and $n^2+n$~inequalities (not counting sign constraints and integer requirements). Use these qualities to reduce the number of variables and obtain a simpler constraint matrix $A$. Then the feasible set is as follows (without sign constraints and integer requirement):
\begin{equation}\label{eq:kmed}
\left[
\begin{matrix}
  &  &                        & &\dots &e_n &e_n& \dots & e_n\\
  &    &   I_{(n-1)n}     &  &\dots &-e_n& 0_n&\dots & 0_n \\
    &    &                   &   &\dots &0_n& -e_n&\dots & 0_n \\  
 -I_n    & -I_n   & \dots& -I_n       &\dots &0_n& 0_n&\dots & -e_n \\       
\end{matrix}
\right]X\leq b,
\end{equation}
where 
\begin{itemize}
\item $b\in\R^{(n^2-1)}$ contains  integral numbers only ($1$, $-1$, $0$ or $k$);
\item $X\in\R^{n^2-1}$ is the vector of decision variables;
\item $e_n=(1,1,\dots,1)^T\in\R^n$;  
\item $I_n$ is an identity matrix of size $n$;
\item $I_{(n-1)n}$ is an identity matrix of size $(n-1)n$;
\item the system matrix $A\in\R^{n^2\times(n^2-1)}$. 
\end{itemize}

Matrix~$A$ contains $n$ blocks of rows ($n$ rows in each block). Add the first $(n-1)$ blocks to the final block of rows (keeping the same order of rows in each block). By doing this, the obtained right hand side vector remains integer, while the last block of rows consists of zeros. If the problem is feasible (that is $k\geq 1$), the final block of rows can be removed. Then the remaining system matrix is 
\begin{equation}\label{eq:B_kmedois}
B=\left[
I_{n(n-1)} C
\right],
\end{equation} 
where 
$$
C=
\left[ 
\begin{matrix}
e_n& e_n&e_n&\dots &e_n&e_n\\
-e_n&0_n&0_n&\dots &0_n&0_n\\
0_n& -e_n&0_n&\dots &0_n&0_n\\
0_n& 0_n&0_n&\dots &-e_n&0_n\\
\end{matrix}
\right],
$$
where $$e_n=(1,1,\dots,1)^T\in\R^n,~0_n=(0,0,\dots,0)^T\in\R^n,~C\in\R^{n(n-1)\times n-1}.$$ 
To complete the proof, it is enough to show that matrix $C$ is totally unimodular.   
Indeed, for any collection of $m$ columns in $C$ ($m\leq n$), the columns can be split into two parts: it is enough to assign any $(m-1)/2$~columns (when $m$ is odd) or $m/2$~columns (when $m$ is even) to one of the parts and the remaining columns in the other part. Then the sum of the columns in part~I minus the sum of the columns in part~II contains $1$, $-1$ or $0$. 

Combine these results with Theorem~\ref{thm:main}. This  completes the proof.

\hskip300pt $\square$

 Therefore, it is enough to apply the simplex method to solve this problem and obtain an integer optimal solution.

 Note that Theorem~\ref{thm:kmedoid} is an important result, since it allows one to avoid integer solvers when applying $k$-medoid method. This results is of interest of cluster analysis and allocation  and therefore has many other potential applications. 
 
 To our best knowledge, there is no result in the literature confirming that all the vertices of the linear relaxations of $k$-medoid formulations are integers. One relevant study~\cite{arhiv:experiments2015} conducts a comprehensive numerical study on $k$-means and $k$-medoid, where most experimental results confirm that the relaxation produce integer (or nearly integer) results: ``LP relaxation remains integral with high probability\rq{}\rq{}. In the same paper,  the authors talk about ``generically unique solutions\rq{}\rq{}, since ``no constraint is parallel to the objective function \rq{}\rq{}. Our approach can deal with the problems where there are such constraints and therefore our approach is more general.  Moreover, the current paper provides an analytical proof that the vertices are integral and therefore the classical implementation of the simplex method always terminates at an integer solution and hence the corresponding problems can be solved efficiently.
 
A comprehensive survey of modern location science~\cite{Location:2015} considers similar problems ($p$-median). 
The authors suggest to solve these problems using a heuristic method. This also confirms that our result is new, since the authors of this recent and comprehensive publication are not aware about the fact that the simplex method obtains an integer solution which is optimal.  
\section{Numerical Experiments}\label{sec:ne}

To demonstrate the applicability of these RLP models, a hypothetical disaster relief problems have been solved using CPLEX Solver 12.7 on a 3.4 GHz processor with 4 GB of RAM.

The results obtained for {expendable and non-expendable} resource allocation were comparable for both implementations: through the simplex method (LP) and branch and bounds (ILP). Namely, the computational time and the optimal value of the cost function were almost the same. This confirms that the branch and bound method terminates after just very few branching stages.

 In the case of the relief centre location experiments, the situation is different, since the corresponding optimisation problems are more complex. Table~\ref{table:experiments}  contains the results of the numerical experiments. In this table we compare  the computational time for the simplex method and IP implementations, the optimal cost value is the same for both implementation and therefore only one column is provided. $k$ represents the number of clusters.  These results indicate that
 the computational time is significantly lower for the simplex method.    
On the top of this, the simplex method reached integer solutions in all the experiments. All these confirm that it is much more efficient to use the classical simplex method in $k$-medoid models.
 \begin{table}
\caption{Relief centre allocation experiments}
 \begin{small}
 \label{table:experiments}
 \begin{tabular}{|c|c|c|c|c|}
 \hline
    Number of data points & $k$&Cost function value&Time (simple)&Time (IP)         \\
    \hline
20&	10&	512.3478602&	0.297&0.562\\
20&	5&	311.3343592&	0.39	&0.391\\
40&	10&	679.4089905&	0.063&	0.078\\
40&	5&	425.939219	&0.063	&0.067\\
60&	10&	848.5919719&	0.14	&0.235\\
60&	5&	490.5917792&	0.187&	0.191\\
80&	10&	961.283263	&0.235	&0.359	\\
80&	5&	565.3424002&	0.25	&0.351	\\
100&	10&	1109.977437&	0.406&	0.594	\\
100&	5&	640.0424524&	0.438&	0.532	\\
120&	10&	1173.702679&	0.594&	0.797	\\
120&	5&	676.7603006&	0.64	&0.718	\\
140&	10&	1241.879626&	0.765&	1.078	\\
140&	5&	715.5979726&	0.86	&1	\\
160&	10&	1297.640491	&1.016	&1.453\\
160&	5&	770.3419758&	1.109&	1.36\\
180&	10&	1346.426155&	1.328&	1.828\\
180&	5&	801.1360661&	1.219&	1.765\\
200&	10&	1425.259571&	1.687&	2.344\\
200&	5&	860.7051955&	1.453&	2.438\\
400&	10&	1954.177312&	8.844&	11.672\\
400&	5&	1170.263786&	7.922&	11.172\\
600&	10&	2379.228282&	25.609&	32.281\\
600&	5&	1406.255502&	20.422&	33.812\\
800&	10&	2729.268892&	48.281&	104.875\\
800&	5&	1617.268652	&43.688&	72.808	\\
\hline
   \end{tabular}
 \end{small}   
            \end{table}
\begin{comment}
It is important to note that most commercial packages implement a number of heuristic techniques in order to accelerate the simplex method. Some of these approaches lead to solutions that are not vertices of the feasible set and therefore are not necessarily integers. We have observed this results in MATLAB when datasets contain more than 100 points. CPLEX implementations of the simplex method are more stable than those implemented in MATLAB and therefore CPLEX was chosen for our numerical experiments.  
\end{comment}



%
%

%

\section{Discussion and future research directions}\label{sec:conclusions}
In this paper we  study two groups of disaster management problems: resource allocation ({expendable and non-expendable}) and relief centre location. We demonstrate that  these problems can be formulated as integer linear programming problems, whose optimal solution can be found by applying the simplex method to the corresponding linear relaxations. There are two main advantage of this approach. Firstly,  most optimisation packages include the simplex method. In particular, small-medium size linear problems can be solved with Excel, while mixed-integer solvers (for example Branch and Bound) may not be available. Secondly, our approach allows one to reduce  the computational time. This is especially prominent in the case of the relief centre location problems, where the corresponding problems are binary.


$k$-medoid method is a very popular clustering technique that is also commonly used in other application. Therefore, another important contribution of this paper is the demonstration that this problem can be solved without application of an integer solver, that is all the computations can be done by applying the classical simplex method, available in most linear packages. This result is achieved by demonstrating that the constraint matrices of the corresponding integer linear programming problems are totally unimodular.

%


One of our future  research direction is to include capacity constraints for the relief centre location problem. Another research direction is to study a modified problem where none of the demand points can be supplied from more than one relief centre. Another future research direction is to extend the results to the case of $p$-median problems.

\section*{Acknowledgement}
This research was supported by the Australian Research Council (ARC),  Solving hard Chebyshev approximation problems through nonsmooth analysis (Discovery Project DP180100602).



\end{document}